\documentclass[a4paper,12pt,twoside,leqno,final]{amsart}
\usepackage{amsmath}
\usepackage{amssymb}

\newtheorem{thm}{Theorem}[section]

\newtheorem{prop}[thm]{Proposition}
\newtheorem{cor}[thm]{Corollary}

\newcommand{\C}{{\mathbb C}}
\newcommand{\D}{{\mathbb D}}
\newcommand{\R}{{\mathbb R}}
\newcommand{\T}{{\mathbb T}}

\newcommand{\eps}{\varepsilon}

\newcommand{\f}{\frac}
\newcommand{\ov}{\overline}

\newcommand{\la}{\lambda}
\newcommand{\ze}{\zeta}
\renewcommand{\th}{\theta}

\newcommand{\ph}{\varphi}

\newcommand{\const}{\text{\rm const}}

\numberwithin{equation}{section}

\title[An extremal problem for Toeplitz operators]
{An extremal problem for functions\\ 
annihilated by a Toeplitz operator}
\author{Konstantin M. Dyakonov}
\address{Departament de Matem\`atiques i Inform\`atica, IMUB, BGSMath, Universitat de Barcelona, Gran Via 585, E-08007 Barcelona, Spain}
\address{ICREA, Pg. Llu\'is Companys 23, E-08010 Barcelona, Spain}
\email{konstantin.dyakonov@icrea.cat}
\thanks{Supported in part by grants MTM2014-51834-P and MTM2017-83499-P from El Ministerio de Econom\'ia y Competitividad (Spain), and by grant 2017-SGR-358 from AGAUR (Generalitat de Catalunya).}

\begin{document}
\begin{abstract}
For a bounded function $\ph$ on the unit circle $\T$, let $T_\ph$ be the associated Toeplitz operator on the Hardy space $H^2$. Assume that the kernel 
$$K_2(\ph):=\{f\in H^2:\,T_\ph f=0\}$$
is nontrivial. Given a unit-norm function $f$ in $K_2(\ph)$, we ask whether an identity of the form $|f|^2=\f12\left(|f_1|^2+|f_2|^2\right)$ may hold a.e. on $\T$ for some $f_1,f_2\in K_2(\ph)$, both of norm $1$ and such that $|f_1|\ne|f_2|$ on a set of positive measure. We then show that such a decomposition is possible if and only if either $f$ or $\ov{z\ph f}$ has a nonconstant inner factor. The proof relies on an intrinsic characterization of the moduli of functions in $K_2(\ph)$, a result which we also extend to $K_p(\ph)$ (the kernel of $T_\ph$ in $H^p$) with $1\le p\le\infty$.
\end{abstract}

\maketitle

\section{Introduction and results}

Let $\T$ stand for the circle $\{\ze\in\C:|\ze|=1\}$ and let $m$ be the normalized Lebesgue measure on $\T$. For $0<p\le\infty$, the space $L^p:=L^p(\T,m)$ will be endowed with the usual norm $\|\cdot\|_p$ (the term \lq\lq quasinorm" should actually be used when $0<p<1$). We also need the {\it Hardy space} $H^p:=H^p(\T)$, viewed as a subspace of $L^p$. The functions in $H^p$ are thus the boundary traces (in the sense of nontangential convergence almost everywhere) of those in $H^p(\D)$, the classical Hardy space on the disk $\D:=\{z\in\C:|z|<1\}$. The latter space consists, by definition, of all holomorphic functions $f$ on $\D$ that satisfy 
$$\sup\left\{\|f_r\|_p:\,0<r<1\right\}<\infty,$$
where $f_r(\ze):=f(r\ze)$ for $\ze\in\T$. 

\par Our starting point is the following observation: Given any $f\in H^2$ with $\|f\|_2=1$, one can find unit-norm functions $f_1,f_2\in H^2$ such that 
\begin{equation}\label{eqn:modfsq}
|f|^2=\f12\left(|f_1|^2+|f_2|^2\right)\,\,\,\text{\rm a.e. on }\T
\end{equation}
and 
\begin{equation}\label{eqn:essdiff}
|f_1|\ne|f_2|\,\,\,\text{\rm on a set of positive measure}. 
\end{equation}
\par This fact (to be explained in a moment) is akin to the classical result that the unit ball of $L^1$ has no extreme points, even though the ball should currently be replaced by a suitable convex subset thereof. Namely, consider the set 
$$V_0:=\left\{|f|^2:\,f\in H^2,\,0<\|f\|_2\le1\right\},$$
i.e., the collection of all functions $g$ on $\T$ that have the form $g=|f|^2$ for some non-null $f$ from the unit ball of $H^2$. We know from basic $H^p$ theory (see \cite[Chapter II]{G}) that the elements $g$ of $V_0$ are characterized by the conditions 
\begin{equation}\label{eqn:glogginlone}
g\ge0\,\,\,\text{\rm a.e. on }\T,\quad g\in L^1,\quad\int_\T\log g\,dm>-\infty,
\end{equation}
and $\|g\|_1\le1$. (The fact that every function $g$ satisfying \eqref{eqn:glogginlone} is writable as $|f|^2$, for some $f\in H^2$, was proved by Szeg\"o in \cite{Sz1}. He then used this representation to study the asymptotic behavior of the polynomials that are orthogonal with respect to such a weight $g$ on $\T$; see \cite[Chapter 12]{Sz2}.) Clearly, the functions $g$ that obey \eqref{eqn:glogginlone} form a convex cone in $L^1$. The portion of that cone lying in the (closed) unit ball of $L^1$ is precisely $V_0$, so this last set is again convex. 
\par We need to show that every function $g\in V_0$ with $\|g\|_1=1$ is a non-extreme point of $V_0$. (By the way, this will imply that $V_0$ has no extreme points at all.) In fact, given such a $g$, we can always find a non-null real-valued function $\tau\in L^\infty$ with the properties that $\int_\T g\tau\,dm=0$ and $\|\tau\|_\infty\le\f12$. This done, we put 
$$g_1:=g(1+\tau),\qquad g_2:=g(1-\tau)$$
and note that 
\begin{equation}\label{eqn:gmidpointonetwo}
g=\f12\left(g_1+g_2\right)\quad\text{\rm a.e. on }\T,
\end{equation}
while $g_1$ and $g_2$ are both in $V_0$. Indeed, for $j=1,2$ we have 
$$\f12g\le g_j\le\f32g\quad\text{\rm a.e. on }\T,$$
which makes \eqref{eqn:glogginlone} true for $g_j$ in place of $g$; also, 
$$\|g_j\|_1=\int_\T g(1\pm\tau)\,dm=\int_\T g\,dm=1.$$ 
Thus \eqref{eqn:gmidpointonetwo} tells us that $g$ is the midpoint of a (nondegenerate) segment whose endpoints $g_j$ ($j=1,2$) lie in $V_0$ and satisfy 
$\|g_j\|_1=1$. Equivalently, a nontrivial decomposition \eqref{eqn:modfsq} with required properties is always possible. 
\par The purpose of this note is to study the equation \eqref{eqn:modfsq} when $f$ lies in a certain subspace of $H^2$, namely, in the kernel of a given Toeplitz operator. The unknowns $f_1$ and $f_2$ are then required to belong to the same subspace and obey \eqref{eqn:essdiff}. Besides, all the functions involved are supposed to be of norm $1$, as before. 
\par For an essentially bounded function $\ph$ on $\T$, we consider the associated operator $T_\ph$ (called the {\it Toeplitz operator with symbol $\ph$}) which acts on $H^2$ by the rule 
$$T_\ph f:=P_+(\ph f),$$
where $P_+$ is the orthogonal projection from $L^2$ onto $H^2$. Assuming that the kernel 
$$K_2(\ph):=\{f\in H^2:\,T_\ph f=0\}$$
is nontrivial, we look at the set 
\begin{equation}\label{eqn:defvphi}
V_\ph:=\left\{|f|^2:\,f\in K_2(\ph),\,0<\|f\|_2\le1\right\},
\end{equation}
i.e., the collection of all functions $g$ on $\T$ that have the form $g=|f|^2$ for some non-null $f$ from the unit ball of $K_2(\ph)$. This set is convex, as we shall soon see, and we are concerned with its (non-)extreme points; once these are determined, we shall arrive at the sought-after information on the solvability of \eqref{eqn:modfsq} for unit-norm functions in $K_2(\ph)$. We are only interested in the case where $\ph\in L^\infty\setminus\{0\}$ (i.e., $\ph$ is non-null), since otherwise $K_2(\ph)=H^2$ and $V_\ph$ reduces to $V_0$, a situation we have already discussed. 
\par Before moving further ahead, we need to gain a better understanding of the moduli of functions in $K_2(\ph)$. This will be achieved by means of Theorem \ref{thm:modulitoep} below, generalizing an earlier result from \cite{DSib}. When dealing with this issue, we temporarily extend our attention to the subspaces 
\begin{equation}\label{eqn:defkerp}
K_p(\ph):=\{f\in H^p:\,T_\ph f=0\}
\end{equation}
with $1\le p\le\infty$, not just with $p=2$. 
\par It should be noted that the operator $P_+$, which kills the function's negative-indexed Fourier coefficients, admits a natural extension to $L^1$, even though 
$P_+(L^1)\not\subset L^1$. This allows us to define the Toeplitz operator $T_\ph$ on $H^1$, whenever $\ph\in L^\infty$, the range of any such operator being contained in $P_+(L^1)$ and hence in every $H^s$ with $0<s<1$. Consequently, the definition \eqref{eqn:defkerp} is meaningful for $p=1$, as well as for all bigger values of $p$. 
\par For a function $f\in H^p$ to be in $K_p(\ph)$, it is necessary and sufficient that the product $\ph f$ be anti-analytic, which in turn amounts to saying that the \lq\lq companion function"
$$\widetilde f:=\ov{z\ph f}$$
is in $H^p$. Thus, 
$$K_p(\ph)=\{f\in H^p:\,\widetilde f\in H^p\},\qquad 1\le p\le\infty.$$
(As a byproduct of this characterization, we mention the fact that every Toeplitz kernel $K_p(\ph)$ enjoys the $F$-property of Havin; see \cite{Hav}.) 
\par A bit more terminology and notation will be needed. Given a nonnegative function $w$ on $\T$ with $\log w\in L^1$, the corresponding {\it outer function} $\mathcal O_w$ is defined a.e. on $\T$ by 
$$\mathcal O_w:=\exp\left\{\log w+i\mathcal H(\log w)\right\},$$
where $\mathcal H$ stands for the harmonic conjugation operator. Functions of the form $\la\mathcal O_w$, with $w$ as above and $\la$ a unimodular constant, will also be referred to as outer. It is well known (and easy to check) that $\mathcal O_w$ extends analytically into $\D$ and has modulus $w$ a.e. on $\T$. Furthermore, $\mathcal O_w\in H^p$ if and only if $w\in L^p$. Finally, we recall that an $H^\infty$ function is said to be {\it inner} if its modulus equals $1$ a.e. on $\T$. See, e.g., \cite[Chapter II]{G} for a systematic treatment of these concepts and of the basic facts related to them.

\begin{thm}\label{thm:modulitoep} Let $1\le p\le\infty$, and suppose $\ph$ is a non-null function in $L^\infty$ for which $K_p(\ph)\ne\{0\}$. Also, let $g$ be a nonnegative function in $L^{p/2}$. The following conditions are then equivalent. 
\par {\rm (i.1)} There is an $f\in K_p(\ph)$ such that $|f|^2=g$ a.e. on $\T$.
\par {\rm (ii.1)} $\ov z\ov\ph g\in H^{p/2}$.
\par\noindent Moreover, if {\rm (ii.1)} holds (with $g$ non-null) and if $I$ is the inner factor of $\ov z\ov\ph g$, then the general form of a function $f\in K_p(\ph)$ with $|f|^2=g$ is given by $f=\mathcal O_{\sqrt g}J$, where $J$ is an inner divisor of $I$. 
\end{thm}

\par Among the possible symbols $\ph$ of our Toeplitz operators, we may single out those which are complex conjugates of inner functions. For such $\ph$'s (i.e., for $\ph=\ov\th$ with $\th$ inner), the corresponding Toeplitz kernels $K_p(\ph)$ take the form 
$$H^p\cap\ov z\th\ov{H^p}=:K^p_\th$$
and are known as {\it star-invariant} or {\it model} subspaces. When $1\le p<\infty$, these are precisely the invariant subspaces of the backward shift operator $f\mapsto(f-f(0))/z$ in $H^p$; see \cite{DSS, N}. We mention in passing that there are deeper connections between the two types of spaces, $K_p(\ph)$ and $K^p_\th$. Namely, in a way, generic Toeplitz kernels can be cooked up from model subspaces; see \cite{DJFA, Hay, Sar} for details.

\par It was in the $K^p_\th$ setting that Theorem \ref{thm:modulitoep} originally appeared in \cite{DSib}; see also \cite[Lemma 5]{DJAM}. In the (tiny) special case where $\ov\ph=\th=z^{n+1}$, the subspace in question is populated by polynomials of degree at most $n$, and the equivalence between (i.1) and (ii.1) above reduces to the classical Fej\'er--Riesz theorem that describes the moduli of such polynomials on $\T$ (see, e.g., \cite[p.\,26]{Sim}). 

\par The role of Theorem \ref{thm:modulitoep} in the present context consists in providing a useful -- and usable -- description of the set $V_\ph$, as defined by \eqref{eqn:defvphi}. Specifically, it tells us that a nonnegative function $g\in L^1\setminus\{0\}$ is in $V_\ph$ if and only if it satisfies $\ov z\ov\ph g\in H^1$ and $\|g\|_1\le1$. This criterion (which implies the convexity of $V_\ph$, among other things) will be repeatedly used hereafter.

\par Our main result, to be stated next, characterizes the extreme points of $V_\ph$. Of course, every function $g\in V_\ph$ with $\|g\|_1<1$ is non-extreme, since 
$$g=\f12(1+\eps)g+\f12(1-\eps)g$$
and $(1\pm\eps)g\in V_\ph$ for suitably small $\eps>0$. Therefore, we only need to consider the case where $\|g\|_1=1$.

\begin{thm}\label{thm:extrpoints} Let $\ph\in L^\infty\setminus\{0\}$ be such that $K_2(\ph)\ne\{0\}$, and let $g\in V_\ph$ be a function with $\|g\|_1=1$. The following are equivalent. 
\par {\rm (i.2)} $g$ is an extreme point of $V_\ph$. 
\par {\rm (ii.2)} $\ov z\ov\ph g$ is an outer function in $H^1$. 
\end{thm}

\par This characterization is reminiscent of de Leeuw and Rudin's theorem (see \cite{dLR} or \cite[Chapter IV]{G}) which identifies the extreme points of the unit ball of $H^1$ as unit-norm outer functions. We also mention a related result from \cite{DPAMS} that describes the extreme points of the unit ball in $K_1(\ph)$. (Namely, these are shown to be the unit-norm functions $f\in K_1(\ph)$ with the property that the inner factors of $f$ and $\widetilde f$ are relatively prime.) In this connection, see also \cite{DKal} and \cite{DMRL}. 
\par We now state a consequence of Theorem \ref{thm:extrpoints} which provides an additional piece of information on the geometry of $V_\ph$. 

\begin{cor}\label{cor:midpoint} Let $\ph\in L^\infty\setminus\{0\}$ be such that $K_2(\ph)\ne\{0\}$. Then every function $g\in V_\ph$ with $\|g\|_1=1$ has the form $g=\f12(g_1+g_2)$, where $g_1$ and $g_2$ are extreme points of $V_\ph$.
\end{cor}

\par Going back to Theorem \ref{thm:extrpoints}, we remark that for $g\in V_\ph$, the function $\ov z\ov\ph g$ can be written as $f\widetilde f$, where $f$ is some (any) function in $K_2(\ph)$ with $|f|^2=g$ and $\widetilde f:=\ov{z\ph f}(\in H^2)$. Consequently, we may rephrase condition (ii.2) above by saying that both $f$ and $\widetilde f$ are outer functions. This in turn leads us to a reformulation of Theorem \ref{thm:extrpoints}, which is perhaps better suited for answering our original question. 

\begin{thm}\label{thm:rephrased} Suppose that $\ph\in L^\infty\setminus\{0\}$, $f\in K_2(\ph)$ and $\|f\|_2=1$. In order that there exist a decomposition of the form \eqref{eqn:modfsq} with some unit-norm functions $f_1,f_2\in K_2(\ph)$ satisfying \eqref{eqn:essdiff}, it is necessary and sufficient that either $f$ or $\widetilde f$ have a nonconstant inner factor.
\end{thm}

\par Finally, we take yet another look at condition (ii.2) of Theorem \ref{thm:extrpoints}. Assuming that (ii.2) holds, we know from Theorem \ref{thm:modulitoep} that the only functions $f\in K_2(\ph)$ with $|f|^2=g$ are constant multiples of $\mathcal O_{\sqrt g}$. It turns out that a similar conclusion is valid, under condition (ii.2), for those functions $f\in K_2(\ph)$ which are merely dominated by $\sqrt g$ in the sense that 
\begin{equation}\label{eqn:domination}
\int_\T\f{|f|}{\sqrt g}dm<\infty.
\end{equation}
(Clearly, this holds in particular when $|f|^2\le\const\cdot g$ on $\T$, let alone when $|f|^2=g$.) In the more general context of $H^p$ spaces with $p\ge1$, the underlying rigidity phenomenon manifests itself in essentially the same way. 

\begin{prop}\label{prop:rigidity} Let $1\le p\le\infty$ and let $\ph\in L^\infty\setminus\{0\}$ be such that $K_p(\ph)\ne\{0\}$. Suppose $g$ is a nonnegative function in $L^{p/2}\setminus\{0\}$ for which $\ov z\ov\ph g$ is an outer function in $H^{p/2}$. Then every function $f\in K_p(\ph)$ satisfying \eqref{eqn:domination} is of the form $f=c\mathcal O_{\sqrt g}$ for some constant $c\in\C$.
\end{prop}

\par In the special case where $\ph$ is the conjugate of an inner function, a similar rigidity result can be found (in a somewhat weaker form) as Theorem 5 in \cite{DSib}. 
\par Now let us turn to the proofs of our current results. 

\section{Proof of Theorem \ref{thm:modulitoep}}

If (i.1) holds, then 
$$\ov z\ov\ph g=\ov z\ov\ph|f|^2=f\cdot\ov z\ov\ph\ov f=f\widetilde f.$$
This last product is in $H^{p/2}$, because $f$ and $\widetilde f$ are both in $H^p$, and we arrive at (ii.1). 
\par Before proceeding to prove the converse, we pause to observe that 
\begin{equation}\label{eqn:logphinlone}
\log|\ph|\in L^1, 
\end{equation}
thanks to our hypotheses on $\ph$. Indeed, let $f_0$ be a non-null function in $K_p(\ph)$. Then $\widetilde f_0:=\ov z\ov\ph\ov f_0$ is in $H^p\setminus\{0\}$ (recall that $|\ph|>0$ on a set of positive measure, while $|f_0|>0$ a.e.), and so 
$$\log|\ph|=\log|\widetilde f_0|-\log|f_0|\in L^1.$$
\par Now suppose that (ii.1) holds, so that $\ov z\ov\ph g=:\mathcal G$ is in $H^{p/2}$. The case of $g\equiv0$ being trivial, we shall henceforth assume that $g$ is non-null; the same is then true for $\mathcal G$. (The latter conclusion relies on \eqref{eqn:logphinlone}, which guarantees that $|\ph|>0$ a.e. on $\T$.) It follows that $\log|\mathcal G|\in L^1$, and hence 
$$\log\,g=\log|\mathcal G|-\log|\ph|\in L^1$$
(where \eqref{eqn:logphinlone} has been used again). We may then consider the outer function $\mathcal O_{\sqrt g}=:F$, so that $F\in H^p$ and $|F|=\sqrt g$, and we further claim that $F\in K_p(\ph)$. 
\par To see why, note that $|\mathcal G|=|\ph|g$, and consequently, the outer factor of $\mathcal G$ equals 
$$\mathcal O_{|\mathcal G|}=\mathcal O_{|\ph|}\mathcal O_g=\Phi F^2,$$
where $\Phi:=\mathcal O_{|\ph|}(\in H^\infty)$. Therefore, letting $I$ denote the inner factor of $\mathcal G$, we have 
\begin{equation}\label{eqn:factg}
\mathcal G=\Phi F^2I.
\end{equation}
Using \eqref{eqn:factg} and the fact that 
\begin{equation}\label{eqn:gfbarf}
g=|F|^2(=\ov FF),
\end{equation}
we now rewrite the identity $\ov z\ov\ph g=\mathcal G$ in the form 
$$\ov z\ov\ph\ov FF=\Phi F^2I,$$
or equivalently, 
\begin{equation}\label{eqn:capftilde}
\ov z\ov\ph\ov F=\Phi FI.
\end{equation}
Thus, $\widetilde F:=\ov z\ov\ph\ov F$ is in $H^p$, which means that $F\in K_p(\ph)$, as claimed above. Finally, we recall \eqref{eqn:gfbarf} to arrive at (i.1), with $f=F$. The equivalence of (i.1) and (ii.1) is thereby verified. 
\par To prove the last assertion of the theorem, assume that $g$ satisfies (ii.1) and that $f$ is an $H^p$ function with $|f|^2=g$. The outer factor of $f$ must then agree with $F$, defined as above, so $f=FJ$ for some inner function $J$. Now, in order that $f$ be in $K_p(\ph)$, it is necessary and sufficient that 
\begin{equation}\label{eqn:littleftilde}
\widetilde f(:=\ov z\ov\ph\ov f)\in H^p.
\end{equation}
On the other hand, multiplying both sides of \eqref{eqn:capftilde} by $\ov J$ yields
$$\widetilde f=\ov z\ov\ph\ov F\ov J=\Phi FI\ov J=\Phi FI/J.$$
It follows that \eqref{eqn:littleftilde} holds if and only if $J$ divides $I$, and the proof is complete.

\section{Proofs of Theorem \ref{thm:extrpoints} and Corollary \ref{cor:midpoint}}

\noindent{\it Proof of Theorem \ref{thm:extrpoints}.} We begin by showing that (i.2) implies (ii.2). Suppose that (ii.2) fails, so that the function $\mathcal G:=\ov z\ov\ph g(\in H^1)$ has a nontrivial inner factor, say $u$. Multiplying $u$ by a suitable unimodular constant, if necessary, we may assume that the number $\int_\T gu\,dm$ is purely imaginary (i.e., belongs to $i\R$). Clearly, $\psi:={\rm Re}\,u$ is then a nonconstant real-valued $L^\infty$ function with $\|\psi\|_\infty\le1$; moreover, 
$$\int_\T g\psi\,dm={\rm Re}\int_\T gu\,dm=0.$$
\par Next, we put 
$$g_1:=g(1+\psi),\qquad g_2:=g(1-\psi)$$ 
and we are going to check that 
\begin{equation}\label{eqn:gjinvphi}
g_j\in V_\ph\quad\text{\rm for}\quad j=1,2.
\end{equation}
Indeed, the above-mentioned properties of $\psi$ imply that 
$$g\left(1\pm\psi\right)\ge0\quad\text{\rm a.e. on}\,\,\T,$$
while 
$$\int_\T g\left(1\pm\psi\right)\,dm=\int_\T g\,dm=1.$$
Thus, $g_1$ and $g_2$ are nonnegative $L^1$ functions, both of norm $1$. We also claim that
\begin{equation}\label{eqn:claimhone}
\ov z\ov\ph g_j\in H^1\quad\text{\rm for}\quad j=1,2.
\end{equation}
To see why, write $G$ for the outer factor of $\mathcal G$ (so that $\mathcal G=Gu$) and note that 
\begin{equation}\label{eqn:goneinhone}
\begin{aligned}
\ov z\ov\ph g_1&=\ov z\ov\ph g(1+\psi)=\mathcal G(1+\psi)\\
&=Gu\left(1+\f12u+\f12\ov u\right)=\f12G(1+u)^2.
\end{aligned}
\end{equation}
A similar calculation yields 
\begin{equation}\label{eqn:gtwoinhone}
\ov z\ov\ph g_2=-\f12G(1-u)^2.
\end{equation}
Because $G\in H^1$ and $(1\pm u)^2\in H^\infty$, the right-hand sides of \eqref{eqn:goneinhone} and \eqref{eqn:gtwoinhone} are both in $H^1$. The claim \eqref{eqn:claimhone} is thereby established, and so is \eqref{eqn:gjinvphi}. 
\par Finally, $g_1\not\equiv g_2$ because $g>0$ a.e. and $\psi$ is non-null. The representation 
\begin{equation}\label{eqn:repronetwo}
g=\f12\left(g_1+g_2\right)
\end{equation}
now allows us to conclude that $g$ is a non-extreme point of $V_\ph$, in contradiction with (i.2). 
\par Conversely, suppose that (ii.2) is fulfilled. Thus, $\mathcal G:=\ov z\ov\ph g$ is an outer function in $H^1$. Now assume that \eqref{eqn:repronetwo} holds with some $g_1$ and $g_2$ in $V_\ph$; hence, in particular, 
\begin{equation}\label{eqn:bothnormsone}
\|g_1\|_1=\|g_2\|_1=1.
\end{equation}
Setting $h:=g_1-g$ and $\mathcal H:=\ov z\ov\ph h$, we further observe that 
\begin{equation}\label{eqn:hinhone}
\mathcal H=\ov z\ov\ph g_1-\ov z\ov\ph g\in H^1.
\end{equation}
Also, we have $g_1=g+h$ and $g_2=g-h$, so \eqref{eqn:bothnormsone} takes the form 
\begin{equation}\label{eqn:gpmhnormone}
\|g+h\|_1=\|g-h\|_1=1.
\end{equation}
Therefore, 
$$\int_\T\left(|g+h|+|g-h|\right)dm=2,$$
or equivalently, 
\begin{equation}\label{eqn:intequalstwo}
\int_\T\left(|1+\Psi|+|1-\Psi|\right)d\mu=2,
\end{equation}
where $\Psi:=h/g$ and $d\mu:=g\,dm$. Because $\mu$ is a probability measure on $\T$ which has the same null-sets as $m$, we may couple \eqref{eqn:intequalstwo} with the obvious inequality 
$$|1+\Psi|+|1-\Psi|\ge2$$
to deduce that we actually have 
$$|1+\Psi|+|1-\Psi|=2$$
a.e. on $\T$. This in turn means that $\Psi$ takes its values in the (real) interval $[-1,1]$. 
\par On the other hand, 
\begin{equation}\label{eqn:psihoverg}
\Psi:=\f hg=\f{\ov z\ov\ph h}{\ov z\ov\ph g}=\f{\mathcal H}{\mathcal G}.
\end{equation}
Recalling that $\mathcal H\in H^1$ (as \eqref{eqn:hinhone} tells us), while $\mathcal G$ is outer, we deduce from \eqref{eqn:psihoverg} that $\Psi$ belongs to the Smirnov class $N^+$ (see \cite[Chapter II]{G}). We also know that $\Psi$ is bounded, whence 
$$\Psi\in N^+\cap L^\infty=H^\infty;$$
and since the only real-valued functions in $H^\infty$ are constants, it follows that $\Psi\equiv c$ for some constant $c\in[-1,1]$. Consequently, $h=cg$ and 
$$\|g\pm h\|_1=(1\pm c)\|g\|_1=1\pm c.$$
Comparing this with \eqref{eqn:gpmhnormone}, we finally conclude that $c=0$. Thus, $h\equiv0$ and $g_1=g_2=g$, so that the only decomposition of the form \eqref{eqn:repronetwo} is the trivial one. This brings us to (i.2) and completes the proof. \qed

\medskip\noindent{\it Proof of Corollary \ref{cor:midpoint}.} If $g$ is an extreme point of $V_\ph$, then it suffices to take $g_1=g_2=g$. Now, if $g$ is non-extreme (so that $\ov z\ov\ph g$ is non-outer), then we may use the representation \eqref{eqn:repronetwo} from the proof of the (i.2)$\implies$(ii.2) part above. To check that the functions $g_1$ and $g_2$ constructed there are actually extreme points of $V_\ph$, we invoke the (ii.2)$\implies$(i.2) part of the theorem, coupled with the fact that the functions $\ov z\ov\ph g_j$ ($j=1,2$) are both outer. The latter is readily seen from \eqref{eqn:goneinhone} and \eqref{eqn:gtwoinhone}, since each of these identities has an outer function, namely $\pm\f12G(1\pm u)^2$, for the right-hand side. \qed

\section{Proof of Proposition \ref{prop:rigidity}}

Let $f\in K_p(\ph)$ be a function satisfying \eqref{eqn:domination}. Then 
\begin{equation}\label{eqn:ratio}
\f{|f|^2}g=\f{\ov z\ov\ph \ov f\cdot f}{\ov z\ov\ph g}=\f{\widetilde ff}{\mathcal G},
\end{equation}
where we write $\widetilde f:=\ov{z\ph f}$ and $\mathcal G:=\ov z\ov\ph g$, as before. Since $\widetilde f$ (as well as $f$) is in $H^p$, while $\mathcal G$ is an {\it outer} function in $H^{p/2}$, it follows that the quotient on the right-hand side of \eqref{eqn:ratio} lies in the Smirnov class $N^+$. The same is therefore true for the left-hand side of \eqref{eqn:ratio}, that is, for $|f|^2/g$. This last ratio also belongs to $L^{1/2}$, as \eqref{eqn:domination} tells us, and so 
$$\f{|f|^2}g\in N^+\cap L^{1/2}=H^{1/2}.$$
Because the only nonnegative $H^{1/2}$ functions are constants (see \cite[p.\,92]{G}), we infer that 
$$|f|^2=\la g$$ 
for some constant $\la\ge0$. 
\par Now, if $\la=0$, then $f\equiv0$ and we are done. Otherwise, since $\mathcal G$ has no inner part, Theorem \ref{thm:modulitoep} (or rather its final assertion, applied with $\la g$ in place of $g$) allows us to conclude that $f$ agrees, up to a constant factor of modulus 1, with the outer function 
$\mathcal O_{\sqrt{\la g}}\left(=\sqrt\la\mathcal O_{\sqrt g}\right)$. The proof is complete.

\end{document}